\documentclass[11pt]{article}
\usepackage[colorlinks=true,breaklinks=true,citecolor=red,linkcolor=blue,urlcolor=red,unicode]{hyperref}
\usepackage{amsmath,amssymb}

\newcommand{\keywords}[1]{\footnote{\hspace{-7mm}
{\bf\scriptsize Keywords:} {\scriptsize \it #1}}}
\newcommand{\address}[1]{\footnote{\hspace{-7mm}
{\bf\scriptsize Authors addresses:} {\bf\scriptsize #1}}}
\makeatletter
\newtheorem{theorem}{Theorem}[section]
\newtheorem{lemma}{Lemma}[section]
\newtheorem{corollary}{Corollary}[section]

\newtheorem{definition}{Definition}[section]
\newtheorem{remark}{\textit{Remark}}[section]
\newtheorem{example}{\textit{Example}}[section]
\numberwithin{equation}{section}

\def\proof{\noindent{\it Proof.}\hspace{3mm}}

\begin{document}
\title{\textbf{Fixed points of asymptotically nonexpansive mappings with center 0 and applications}}
\author{\textsc{Abdelkader Dehici$^{1}$, Sami Atailia$^{2}$ and Najeh Redjel$^{1}$}}

\date{}
\maketitle
\begin{abstract} In this paper, we investigate the existence of fixed points for asymptotically nonexpansive mappings with center 0 defined on closed convex subsets of various Banach spaces. Three applications are given. Firstly, we prove that our results refine those concerning alternate convexically nonexpansive (in short; ACN) mappings studied by P. N. Dowling in " On a fixed point result of Amini-Harandi in strictly convex Banach spaces, Acta. Math. Hungar., 112 (1-2), (2006), 85-88" . Secondly, by using Lau's result in " Closed convex invariant subsets of $L_p(G)$, Trans. Amer. Math. Soc., {\bf 232}, (1977), 131-142", we give another characterization for the noncompactness of locally compact groups $G$. Finally, we discuss the existence of a solution for a nonlinear transport equation without using compactness results.
\end{abstract}


\address{(1) Department of Mathematics and Computer Science\\
University of Souk-Ahras,P.O.Box 1553, Souk-Ahras 41000, Algeria\\
(2) Department of Mathematics, Faculty of Science, Boumerdes University,\\
Boumerdes 35000, Algeria. \\
\texttt{E-mails: dehicikader@yahoo.fr, s.atailia@univ-boumerdes.dz, najehredjel@yahoo.fr}}
\vskip 0.5 cm
\section {Introduction}
In 1922, S. Banach \cite{Ban} established his metric fixed point theorem, named after him, for contraction mappings. In the case where the contraction constant equals 1 (the nonexpansive case), it is easy to construct examples that do not have fixed points. In fact, it suffices to take the complete  metric space $(X, \widetilde{d})$ where $X = \{0, 1\}, \widetilde{d}$ is the discrete metric and $T: X\longrightarrow X$ is a self-mapping defined by $T(0) = 1$ and $T(1) = 0$. Then $T$ satisfies $\widetilde{d}(Tx, Ty) = \widetilde{d}(x, y)$ for all $x, y \in X$. But $T$ does not have fixed points. So, since the situation does not work in the case of complete metric spaces, what about the framework of Banach spaces?
\vskip 0.3 cm
\noindent In 1965, F. E. Browder, D. G\"{o}hde and W. A. Kirk (\cite{Bro, Goh, Kir}) showed  that every nonexpansive mapping defined on a bounded closed convex subset of a uniformly convex Banach space (or more generally reflexive Banach space having normal structure) has at least a fixed point. Their works were the foundation of the fixed point theory for nonexpansive mappings which illustrated the primary role of the geometry of Banach spaces in this axis of research. It was the birth of an interesting domain of nonlinear functional analysis which attracted the attention of many mathematician. For more details, see \cite{Aks, Bel, Goe2, Goe}.

\vskip 0.3 cm
\noindent A Banach space $X$ is said to have the fixed point property (resp. the weak fixed point property) for nonexpansive mappings (in short; FPP) (resp. in short; w-FPP) if for all nonempty bounded closed (resp. weakly compact) convex subset $C$ of $X$, every nonexpansive mapping $T: C \longrightarrow C$ has at least a fixed point. According to this terminology, we can also define the fixed point property (resp. weak fixed point property) for generalized nonexpansive mappings (see \cite{Bae, Bet, Bet1, Bog, Fal1, Goe2, Rei3, Rei4, Smy, Suz}).

\vskip 0.3 cm
\noindent One of the most passionate subjects is the link between the FPP and the FPP for generalized nonexpansive mappings. In 1976, C. S. Wong \cite{Won1} proved that the w-FPP for Kannan mappings (mappings that satisfy $\|Tx - Ty\| \leq \displaystyle \frac{1}{2} (\|Tx - x\| + \|Ty - y\|), $ for all $x, y \in C$, see also \cite{Rei1}) characterizing the quasi-weak normal structure which is possessed by strictly convex and separable Banach spaces (see \cite{Won2}). Consequently, Banach space $L^{1}([0, 1])$ has the w-FPP for Kannan mappings. However, according to the famous result due to D. Alspach (see \cite{Als}), $L^{1}([0, 1])$ does not have the w-FPP (resp. FPP). In the same direction, K. K. Tan \cite{Tan} constructed a separable and reflexive Banach space that have FPP but fails to have the FPP for Kannan mappings.

\vskip 0.3 cm
\noindent In general, the link between the w-FPP and the w-FPP for generalized nonexpansive mappings is not known yet where $X$ does not have neither the weak normal structure nor the quasi-weak normal structure. We have just a few contributions in this direction. For this, we cite for example, the works (\cite{Ata, Bet, Bet1, Fal1, Goe1, Fus, Smy, Suz}). We remind that in the case of real Hilbert spaces, the FPP property characterizes the boundedness of closed convex subsets. This result was proved by W. O. Ray \cite{Ray} and then simplified by R. Sine \cite{Sin}. Also, if $X = c_0$ (the Banach space of real sequences that converge to zero)., the FPP characterizes the weak compactness of bounded closed convex subsets. Recently, the authors in \cite{Ata, Deh, Deh1} gave an investigation on the FPP for $(c)$-mappings in the bounded and unbounded cases. Particularly, by using the results of W. Takahashi et al \cite{Tak}. A. Dehici and S. Atailia \cite{Deh} proved a variant result of Ray for $(c)$-mappings and the problem in the case of abstract Banach spaces, is still open.

\vskip 0.3 cm
\noindent In \cite{Fal2}, J. Garc\'{\i}a-Falset et al introduced the class of nonexpansive mappings with center as an extension of quasi-nonexpansive mappings. They proved by examples that a center is not necessarily a fixed point for these considered mappings. The advantage of the contributions in \cite{Fal1} concerning the study of fixed points, is that they are established for mappings which are not necessarily self-mappings but they are defined on nonempty subsets of various Banach spaces.
\vskip 0.3 cm
\noindent In this paper, we are interested in the class of asymptotically nonexpansive mappings that have zero as a center which contains the class of nonexpansive mappings with center 0. We study the existence of fixed points for these mappings defined on bounded closed convex subsets of reflexive strictly convex Banach spaces and spaces having Kadec-Klee property. Many illustrative examples are given. Afterwards, we materialize our results by three concrete applications.
\vskip 0.3 cm
\noindent For the first one, we refine the results of A. Amini-Harandi \cite{Ami} and P. Dowling \cite{Dow} which are established in the case of weakly compact convex subsets of strictly convex Banach spaces and we extend them to the case of closed convex (not necessarily bounded) subsets of reflexive strictly convex Banach spaces.
\vskip 0.3 cm
\noindent For the second one, with a result due to A. T-M. Lau \cite{Lau} concerning the characterization of closed convex subsets invariant by modular isometries which are defined on $L_p(G) (1 < p < \infty)$ where $G$ is a locally compact noncompact group, we give a characterization for the noncompactness of $G$ using orbits associated with these isometries.
\vskip 0.3 cm
\noindent The last application is devoted to study the existence of solutions for a nonlinear transport equation with contractive boundary conditions. More precisely, we refine the results of \cite{Lat} in this sense and we show following our assumptions that our results are independent of any use of the compactness argument established in the first section of \cite{Lat}.
\vskip 0.3 cm
\section {Preliminaries and Preparatory Results}
In this paper, we introduce a large class of asymptotically nonexpansive mappings with center 0 containing in particular that of nonexpansive mappings $T: C \longrightarrow C$ having 0 as a fixed point (if $0 \in C$).
\begin{definition}\label{def 2.1}\rm Let $C$ be a nonempty subset of a Banach space $(X, \|.\|)$ and let $T: C \longrightarrow C$ be a self-fmapping. $T$ is said to be nonexpansive if
\begin{equation}
\|Tx- Ty\|\leq \|x- y\| \ \  \hbox{for all} \ \ x,y \in C.
\end{equation}
\end{definition}

\begin{definition}\label{def 2.2}\rm Let $C$ be a nonempty subset of a Banach space $(X, \|.\|)$ and let $T: C \longrightarrow C$ be a self-mapping. $T$ is said to be asymptotically nonexpansive if
\begin{equation}
\displaystyle \limsup_{n \longrightarrow +\infty}\|T^{n}x- T^{n}y\|\leq \|x-y\| \ \ \hbox{for all} \ \ x, y \in C.
\end{equation}
\end{definition}

\begin{remark}\label{rem 2.1}\rm It is easy to see that every nonexpansive mappings is asymptotically nonexpansive while the converse is not true in general as the following example shows:
\end{remark}

\begin{example}\label{exa 2.1}\rm Let $T: [0,1]\longrightarrow [0,1]$ be defined by $Tx= \sqrt{x}$ if $x > 0$ and $T(0)= 1$. $T$ is asymptotically nonexpansive while $T$ is not nonexpansive since $T$ is not continuous at $x_{0}= 0$.
\end{example}

\begin{remark}\label{rem 2.2}\rm As it is indicated in \cite{Lau1}, the notion of asymptotically nonexpansive given in Definition 2.2 is different and more general than the notion of nonexpansive mappings introduced by K. Goebel and W. A. Kirk in \cite{Goe1}.
\end{remark}

\begin{definition}\label{def 2.3}\rm Let $C$ be a nonempty subset of a Banach space $(X, \|.\|)$ and let $T: C \longrightarrow C$ be a self-mapping. $T$ is said to be nonexpansive mapping with center $0$ if $\|Tx\|\leq \|x\|$ for all $x \in C$.
\end{definition}

\begin{remark}\label{rem 2.3}\rm Clearly if $0 \in C$ and $T$ is a nonexpansive mapping having $0$ as a fixed point, then $T$ is nonexpansive with center 0 but the next examples show that there exist nonexpansive mappings with center $0$ which fail to be nonexpansive.
\end{remark}

\begin{example}\label{exa 2.2}\rm Let \begin{align}
T: [0,3] &\longrightarrow [0,3] \nonumber \\
x & \longrightarrow \left\{
                      \begin{array}{ll}
                        0 & \hbox{if} \ x \neq 3;\\
                        1 & \hbox{if} \ x= 3.
                      \end{array}
                    \right. \nonumber
\end{align}

\noindent Obviously, we have $|Tx|\leq |x|$ for all $x \in [0, 3]$. However, $T$ is not nonexpansive since $T$ is not continuous at $x_{0}= 3$.
\end{example}

\begin{example}\label{exa 2.3}\rm Let \begin{align}
T: [0,1] & \longrightarrow [0,1] \nonumber  \\
x & \longrightarrow x^{2} \nonumber
\end{align}
\noindent If $x \in [0,1]$ then $|Tx|\leq |x|$, however $T$ is not nonexpansive. To see this, it suffices to take $x_{1}= \displaystyle \frac{1}{2}$ and $x_{2}= \displaystyle \frac{2}{3}$. For these values, we have
\vskip 0.3 cm
\centerline{$|Tx_{1}- Tx_{2}|= |x_{1}^{2}- x_{2}^{2}|= \left|(\displaystyle \frac{1}{2})^{2}- (\displaystyle \frac{2}{3})^{2}\right|= \left|\displaystyle \frac{1}{2}+ \displaystyle \frac{2}{3}\right|\left|\displaystyle \frac{1}{2}- \displaystyle \frac{2}{3}\right|> |x_{1}- x_{2}|$.}
\vskip 0.3 cm
\end{example}

\begin{remark}\label{rem 2.4}\rm From Remark \ref{rem 2.1}, we observe that nonexpansive mappings with center 0 are asymptotically nonexpansive with the same center.
\end{remark}

\vskip 0.3 cm
\begin{definition}\label{def 2.4}\rm Let $C$ be a nonempty subset of a Banach space $X$ and let $T: C \longrightarrow C$ be a self-mapping. $T$ is said to be a $(c)$-mapping if there exist $a, c \in [0,1] $ with $ (c > 0)$ and $a+ 2c= 1$ such that
\begin{equation}
\|Tx-Ty\|\leq a \|x-y\|+ c(\|x-Ty\|+ \|Tx-y\|) \  \  \hbox{for all} \ \ x,y \in C.
\end{equation}
\end{definition}

\begin{remark}\label{rem 2.5}\rm A simple calculation shows that Example \ref{exa 2.2} is a (c)-mapping for $\displaystyle \frac{1}{3}$. So, $(c)$-mappings can be discontinuous.
\end{remark}

\begin{remark}\label{rem 2.6}\rm It is worth noting that there exit examples which are nonexpansive and $(c)$-mappings at the same time. To see this, it suffices to take $(X, \|.\|)= (\mathbb{R}, |.|)$ and $T: \mathbb{R} \longrightarrow \mathbb{R}$ defined by $Tx= x+ a$ with $a \neq 0$.
\end{remark}

\noindent The following lemma due to the J. S. Bae \cite{Bae} is a useful tool in the investigation of fixed points for $(c)$-mappings.

\begin{lemma}\label{lem 2.1}\rm Let $C$ be a bounded subset of a Banach space $X$ and let $T: C \longrightarrow C$ be a $(c)$-mapping. Then $T$ is asymptotically regular i.e.,
\vskip 0.3 cm
\centerline{$\displaystyle \lim_{n \longrightarrow +\infty}\|T^{n+1}x- T^{n}x\|= 0 \hskip 0.5 cm \hbox{for all} \ x \in C.$}
\vskip 0.3 cm
\end{lemma}

\begin{example}\label{exa 2.4}\rm Let $C= [0,1]\subseteq \mathbb{R}$ where $\mathbb{R}$ is equipped with its usual norm and let $T: [0,1] \longrightarrow [0,1]$ be defined by $Tx= 1-x$. It is easy to see that $T$ is nonexpansive. However, $T$ cannot be a $(c)$-mapping since $T^{2k+1}(0)= 1$ and $T^{2k}(0)= 0$ and the claim follows from Lemma \ref{lem 2.1}.
\end{example}

\begin{definition}\label{def 2.5}\rm Let $C$ be a weakly compact convex subset of a Banach space $X$. $C$ is said to have $(c)$-FPP if every $(c)$-self-mapping on $C$ has a fixed point.
\end{definition}

\begin{remark}\label{rem 2.7}\rm It is an open problem whether $(c)$-FPP holds if FPP is satisfied.
\end{remark}

\begin{definition}\label{def 2.6}\rm A Banach space $X$ is said to be uniformly convex if for each $\epsilon \in (0, 2]$ there exists $\delta > 0$ such that for all $x, y \in X, \|x\| \leq 1, \|y\| \leq 1, \|x - y \| > \epsilon \Longrightarrow \|\displaystyle \frac{x + y}{2}\| \leq \delta$.
\end{definition}

\begin{definition}\label{def 2.7}\rm A Banach space $X$ is said to be strictly convex if for all $x,y \in X, x \neq y$, we have
\vskip 0.3 cm
\centerline{$\|x\|= \|y\|= 1 \Longrightarrow \|\displaystyle \frac{x+ y}{2}\| < 1.$}
\vskip 0.3 cm
\end{definition}

\begin{definition}\label{def 2.8}\rm Every uniformly convex Banach space $X$ is strictly convex while the converse is not true in general. Recall that the Lebesgue spaces $L_p(\mu)$ are uniformly convex.
\end{definition}

\begin{example}\label{exa 2.5}\rm Let $X= C([0,1])$ be the Banach space of scalar continuous functions defined on $[0,1]$ equipped with the following norm
\vskip 0.3 cm
\centerline{$ |\|f\||= \displaystyle \sup_{t \in [0,1]}|f(t)|+ \|f\|_{L^{2}([0,1])}.$}
\vskip 0.3 cm
\noindent Then $\left(X, |\|.\||\right)$ is strictly convex but not uniformly convex Banach space (for more details, see the last paragraph in page 24 of \cite{Goe}).
\end{example}

\begin{definition}\label{def 2.9}\rm Let $C$ be a subset of a Banach space $X$ and let $T: C \longrightarrow C$ be a self-mapping. $T$ is called generalized nonexpansive if there exists $a, b, c \in [0,1]$ such that $a+ 2b+ 2c= 1$ and
\vskip 0.3 cm
\centerline{$\|Tx- Ty\|\leq a \|x-y\|+ b (\|Tx- x\|+ \|Ty - y\|)+ c(\|Tx- y\|+ \|Ty- x\|)$}
\vskip 0.3 cm
\noindent for all $x,y \in C$.
\end{definition}

\begin{remark}\label{rem 2.8}\rm By a simple calculation, we infer that if $0 \in C$ and $T$ is a generalized nonexpansive mapping having $0$ as a fixed point, then necessarily $T$ is nonexpansive with center 0.
\end{remark}

\begin{definition}\label{def 2.10}\rm Let $C$ be a nonempty subset of a Banach space $X$ and let $T: C \longrightarrow C$ be a self-mapping, $T$ is said to be a Suzuki mapping if for all $x, y \in C$
\vskip 0.3 cm
\centerline{$\displaystyle \frac{1}{2}\|Tx- x\|\leq \|x-y\| \Longrightarrow \|Tx- Ty\| \leq \|x- y\|.$}
\vskip 0.3 cm
\end{definition}
\vskip 0.3 cm
\noindent The class of Suzuki mappings was introduced and studied in 2008 by T. Suzuki \cite{Suz} who proved that this class contains strictly the set of nonexpansive mappings. Other interesting properties related to this class can be found in \cite{Fus, Suz}.

\begin{remark}\label{rem 2.9}\rm Similar to Remark \ref{rem 2.8}, if $0 \in C$ and $T$ is a Suzuki mapping having $0$ as a fixed point then $\|Tx\|\leq \|x\|$ for all $x \in C$.
\end{remark}

\begin{definition}\label{def 2.11}\rm (see Definition 2 in \cite{Fal1} and Definition (2) in \cite{Fus}) Let $C$ be a nonempty subset of a Banach space $X$. $T$ is said to satisfy the property $(E_{\mu_0})$ if there exists $\mu_{0} \geq 1$ such that
\vskip 0.3 cm
\centerline{$\|x- Ty\|\leq \mu_{0}\|Tx-x\|+ \|x- y\|$}
\vskip 0.3 cm
\noindent for all $x,y \in C$.
\end{definition}
\vskip 0.3 cm
\noindent Let us give now the following important useful lemma.

\begin{lemma}\label{lem 2.2}\rm (see Lemma 7 in \cite{Suz} and Proposition 3.6 in \cite{Fus}) Let $C$ be a nonempty subset of a Banach space $X$. If $T: C \longrightarrow C$ is a Suzuki mapping or a generalized nonexpansive mapping. Then $T$ satisfies the property $(E_{\mu_{0}})$.
\end{lemma}

\begin{example}\label{exa 2.6}\rm Nonexpansive mappings with center $0$ are not necessarily generalized nonexpansive or Suzuki mappings. Indeed, to see this, take $T_{|[0, \frac{2}{3}]}: [0, \frac{2}{3}]\longrightarrow [0, \frac{2}{3}]$ defined by $Tx = x^{2}$. Then $T$ is nonexpansive with center 0 (see Example 2.3) but it does not satisfy the property $(E_{\mu_{0}})$ (see Example 3.7 in \cite{Fus}).  So, by Lemma 2.2, $T$ is neither a generalized nonexpansive mapping nor a Suzuki mapping.
\end{example}

\section{Main results}
\noindent We start this section by the following main result.

\begin{theorem}\label{th 3.1} (see also page 1207 in \cite{Lau1}) \rm Let $C$ be a nonempty closed convex (not necessarily bounded) subset of a reflexive strictly convex Banach space. Then, there exists a unique element $z_0 \in C$ such that
\vskip 0.3 cm
\centerline{$ \|z_0\| = \delta_0= \inf \{\|y\|: y \in C\}.$}
\end{theorem}
\proof {\it First case:} $0 \in C$, then $z_0 = 0$.
\vskip 0.3 cm
\noindent {\it Second case:} $0 \notin C$.  Let us define the following set
\vskip 0.3 cm
\centerline{$ C_0 = \{z \in C: \|z\| = \delta_0\}$.}
\vskip 0.3 cm
\noindent It is easy to see that $C_0 = \displaystyle \bigcap_{\epsilon > 0} C_{0, \epsilon}(\delta_0)$ where $C_{0, \epsilon}(\delta_0) = \{z \in C: \|z\| \leq \delta_0 + \epsilon\}$. Since $X$ is reflexive then for all $\epsilon > 0, C_{0, \epsilon}$ is a nonempty weakly compact convex subset of $X$ and by the finite intersection property, $C_0$ is a nonempty weakly compact convex subset of $X$.
\vskip 0.3 cm
\noindent Now, let $z_1, z_2 \in C_0$, with $z_1 \neq z_2$. Since $X$ is strictly convex, then necessarily we obtain that
\vskip 0.3 cm
\centerline{$\frac{\|z_1 + z_2\|}{2} < \frac{1}{2} (\|z_1\| + \|z_2\|) = \delta_0$.}
\vskip 0.3 cm
\noindent which is a contradiction. Consequently, the set $C_0$ is reduced to a singleton $\{z_0\}$ which is the desired result.

\vskip 0.3 cm
\begin{remark}\label{re 3.1}\rm In the above theorem, the element $z_0$ is called the metric projection of 0 on $C$ and is denoted by $P_C(0)$. In the case where $X$ is a real Hilbert space then $P_C(0)$ is characterized by the unique element $z_0 \in C$ such that
\vskip 0.3 cm
\centerline{$Re <z_0, z - z_0> \geq 0, z \in C.$}
\end{remark}
\vskip 0.3 cm
\noindent From the proof of Theorem 3.1, we can deduce the following

\vskip 0.3 cm
\begin{corollary}\label{co 3.1}\rm Let $C$ be a weakly compact convex subset of a strictly convex Banach space (not necessarily reflexive). Then, there exists a unique element $\widehat{z_0} \in C$ such that
\vskip 0.3 cm
\centerline{$ \|\widehat{z_0}\| = \widetilde{\delta_0}= \inf \{\|y\|: y \in C\}.$}
\end{corollary}
\vskip 0.3 cm
\noindent The first fixed point result concerning asymptotically nonexpansive mappings with center 0 is given in the following.
\vskip 0.3 cm
\begin{theorem}\label{th 3.1}\rm Let $C$ be a closed convex subset of a reflexive strictly convex Banach space and let $T: C \longrightarrow C$ be an asymptotically nonexpansive mapping with center $0$.
\begin{description}
  \item[$(\imath)$] If $0 \in C$ and $T$ is continuous then $0$ is a fixed point for $T$.
  \item[$(\imath\imath)$] If $0 \notin C$ and $T$ is weakly continuous then $T$ has a fixed point in $C$.
\end{description}
\end{theorem}

\proof Since $T$ is asymptotically nonexpansive, then
\vskip 0.3 cm
\centerline{$\displaystyle \limsup_n \|T^{n}x\| \leq \|x\| \ \ \hbox{for all} \ \ x \in C.$}
\vskip 0.3 cm
\begin{description}
  \item[$(\imath)$] If $0 \in C$, we have
\vskip 0.3 cm
\centerline{$\displaystyle \limsup_n \|T^{n}0\|\leq 0,$}
 \vskip 0.3 cm
\noindent which leads to
\vskip 0.3 cm
\centerline{$\displaystyle \lim_{n \longrightarrow + \infty}  \|T^{n}0\|= 0.$}
 \vskip 0.3 cm
\noindent So,
\vskip 0.3 cm
\centerline{$\displaystyle \lim_{n \longrightarrow + \infty}T^{n}0= 0,$}
\vskip 0.3 cm
\noindent since $T$ is continuous, we get
\vskip 0.3 cm
\centerline{$\displaystyle \lim_{n \longrightarrow + \infty}T^{n+ 1}0= T(0),$}
\vskip 0.3 cm
\noindent consequently, we have $T(0)= 0$.
  \item[$(\imath\imath)$] If $0 \notin C$, following Theorem 3.1 there exists a unique $x_{0}\in C$ such that
\vskip 0.3 cm
\centerline{$\|x_{0}\|= \inf\{\|x\|: x \in C\}.$}
\vskip 0.3 cm
So
\vskip 0.3 cm
\centerline{$\displaystyle \limsup_n \|T^{n}x_{0}\|\leq \|x_{0}\|.$}
\vskip 0.3 cm
\noindent Assume that $y_0 \in \overline{(T^{n}x_0)_n}^{w}$ (the weak closure of the sequence $(T^{n}x_0)_n$). This implies the existence of a subsequence $(T^{n_k}x_{0})_k$ such that $T^{n_k}x_{0}\longrightarrow y_{0}$ weakly, using the lower semicontinuity of the norm, it follows that
\vskip 0.3 cm
\centerline{$\|y_{0}\|\leq \displaystyle \liminf_k \|T^{n_k}x_{0}\|\leq \displaystyle \limsup_k \|T^{n_k}x_{0}\| \leq \displaystyle \limsup_n \|T^{n}x_{0}\| \leq \|x_{0}\|.$}
\vskip 0.3 cm
\noindent So, necessarily
\vskip 0.3 cm
\centerline{$y_0 \in \{x \in C: \|x\|= \|x_{0}\|\}.$}
\vskip 0.3 cm
\noindent Hence, we deduce that $y_0= x_{0}$. Thus every weakly convergent subsequence $(T^{m_k}x_0)_k$ in $C$ weakly converges to $x_{0}$. Therefore the sequence $(T^{n}x_0)_n$ weakly converges to $x_{0}$. Since $T$ is weakly continuous, we get
\vskip 0.3 cm
\centerline{$T(x_{0})= T( wk-\displaystyle \lim_{n}T^{n}x_{0})= wk-\displaystyle \lim_{n}T^{n+1}x_{0}= x_{0}$}
\end{description}
\noindent (where $wk-\displaystyle \lim_{n}$ is the weak limit) which is the desired result.
\vskip 0.3 cm

\begin{remark}\label{rem 3.2}\rm Following its proof, it can be seen that the assertion 1 of Theorem 3.2 holds  for a nonempty closed convex subset $C$ of an arbitrary Banach space $X$.
\end{remark}
\vskip 0.3 cm
\noindent Let us state the following theorem due to P. Dowling \cite{Dow} (see also \cite{Ami}).

\begin{theorem}\label{th 3.3}\rm Let $C$ be weakly compact convex subset of a strictly convex Banach space and let $T: C \longrightarrow C$ be a nonexpansive mapping with center $0$. Then $T$ has a fixed point.
\end{theorem}
\vskip 0.3 cm
\noindent By using Theorem 3.2 and adapting the same techniques given in the proof of the above theorem (for more details, see Theorem 3 in \cite{Dow}), we can derive the following
\vskip 0.3 cm
\begin{corollary}\label{cor 3.2}\rm Let $C$ be a closed convex subset of a reflexive strictly convex Banach space and let $T: C \longrightarrow C$ be a nonexpansive mapping with center $0$. Then $T$ has a fixed point.
\end{corollary}
\vskip 0.3 cm

\begin{remark}\label{rem 3.3}\rm We observe that Theorem \ref{th 3.1} is more general than Theorem \ref{th 3.3} in the case of reflexive strictly convex spaces since the boundedness of $C$ in Corollary 3.2 is dropped.
\end{remark}

\begin{corollary}\label{cor 3.3}\rm Let $C$ be a closed convex subset of a reflexive strictly convex Banach space $X$. Assume that $T: C \longrightarrow C$ and $S: C \longrightarrow C$ are two self-mappings such that $S$ is into and $\|TS(x)\|\leq \|S(x)\|$ for all $ x \in C$. Then $T$ has a fixed point.
\end{corollary}

\proof The assumption that $\|TS(x)\|\leq \|S(x)\|$ for all $x \in C$ is equivalent to the fact that $T$ is a nonexpansive mapping with center $0$. So the result follows immediately from Corollary 3.2.

\begin{corollary}\label{cor 3.4}\rm Let $C$ be a closed convex subset of a reflexive strictly convex Banach space $X$ and let $T: C \longrightarrow C$ be a self-mapping satisfying that
\vskip 0.3 cm
\hfill{$\|T(x) + T(y)\|\leq \|x + y\|$ for all  $x, y \in C.$ \hfill{} (3.1)}
\vskip 0.3 cm
\noindent Then $T$ has a fixed point in $C$.
\end{corollary}
\proof By taking $x = y$, we observe that $T$ is necessarily nonexpansive with center 0. Now, the result is an immediate consequence of Corollary 3.2.
\vskip 0.3 cm
\noindent In the same way of Corollary 3.4 and using Theorem 3.3, we can derive the following.
\vskip 0.3 cm
\begin{corollary}\label{cor 3.5}\rm Let $C$ be a weakly compact convex subset of a strictly convex Banach space $X$ and let $T: C \longrightarrow C$ be a self-mapping satisfying that
\vskip 0.3 cm
\hfill{$\|T(x) + T(y)\|\leq \|x + y\|,$ for all $x, y \in C.$ \hfill{} (3.2)}
\vskip 0.3 cm
\noindent Then $T$ has a fixed point in $C$.
\end{corollary}

\begin{example}\label{exa 3.1}\rm Let $X = {\mathbb{R}}^{2}$ equipped with the Euclidean norm and let
\vskip 0.3 cm
\centerline{$C = \{(x_1, x_2) \in  {\mathbb{R}}^{2}, x_1, x_2 \geq 0, x_1^{2} + x_2^{2}\leq 1\}.$}
\vskip 0.3 cm
\noindent Define $T: C \longrightarrow C$ by $T(x, y) = (x^{2}, y^{2})$. So, obviously $X$ is strictly convex (since $X$ is uniformly convex) and $C$ is a bounded closed convex subset of $X$. Furthermore, for all $(x_1, y_1), (x_2, y_2) \in C$, we have
\vskip 0.3 cm
\centerline{$\|T(x_1, y_1) + T(x_2, y_2)\| = \sqrt{(x_1^{2} + x_2^{2})^{2} + (y_1^{2} + y_2^{2})^{2}}$}
\vskip 0.3 cm
\centerline{$\leq \sqrt{(x_1 + x_2)^{2} + (y_1 + y_2)^{2}}$}
\vskip 0.3 cm
\centerline{$= \|(x_1, y_1) + (x_2, y_2)\|$}
\vskip 0.3 cm
\noindent Now, all assumptions of Corollary 3.5 are satisfied and the existence of a fixed point for $T$ is ensured. Clearly,  $(0, 0), (1, 0)$ and $(0, 1)$ are the fixed points of $T$ in $C$.
\vskip 0.3 cm
\end{example}

\vskip 0.3 cm
\begin{corollary}\label{cor 3.6}\rm Let $C$ be a weakly compact convex subset of a strictly convex Banach space $X$ and let $T: C \longrightarrow C$ be a self-mapping such that there exists $k_{0}> 1$ for which $T^{k_0}$ is nonexpansive with center $0$. If $T$ is a $(c)$-mapping then $T$ has a fixed point in $C$.
\end{corollary}

\proof Denote $A= T^{k_{0}}$. So by assumption, $A$ is nonexpansive with center $0$. Then, by Theorem 3.3, $A$ has a fixed point $x_{0}\in C$. Then, since $T$ is a $(c)$-mapping then $x_{0}$ is a fixed point in $C$ (see Proposition 4.1 in \cite{Deh}).

\begin{example}\label{exa 3.2}\rm The mapping $T$ in Example 2.2 is a $(c)$-mapping and it is also nonexpansive with center $0$. Obviously, we have $T^{2}= 0$. Then, $T^{2}$ has $0$ as a unique fixed point in $[0,3]$. But, it is easily seen that $0$ is also the unique fixed point of $T$ in $[0,3]$.
\end{example}

\begin{definition}\label{def 3.1}\rm A Banach space $X$ is said to have Kadec-Klee property (in short, KK property) if for every sequence $(x_{n})_{n}$ in $X$ such that if $x_{n}$ converges weakly to $x$ and $\|x_{n}\| \longrightarrow \|x\|$ then $x_{n}$ converges in norm to $x$.
\end{definition}

\begin{remark}\label{rem 3.4}\rm It is easy to deduce that Kadec property (for which weak topology and norm topology are the same) implies KK property but the converse is not true in general (see \cite{Tro}). In particular, spaces having Schur property satisfy KK property.
\end{remark}

\begin{example}\label{exa 3.3}\rm Banach spaces $L_{p}(\mu) \ (1 < p < \infty)$ have KK property.
\end{example}

\begin{remark}\label{rem 3.5}\rm Recall that the classes of strictly convex Banach spaces and those having KK property are different. Indeed $l^{1}(\mathbb{N})$ has Kadec-Klee property but $l^{1}(\mathbb{N})$ is not strictly convex. In addition, the space $(c_{0}, |\|.\||)$ where $|\|.\||$ is defined by
\vskip 0.3 cm
\centerline{$|\|x\||= \|x\|_{\infty}+ (\displaystyle \sum_{n= 1}^{\infty}(\displaystyle \frac{x_{n}}{n})^{2})^{\frac{1}{2}} \ \ \hbox{for all} \ \ x\in c_{0}$}
\vskip 0.3 cm
\noindent is strictly convex but fails to have KK property (see Example 23 in \cite{Fal2}).
\end{remark}

\begin{definition}\label{def 3.2}\rm Let $C$ be a nonempty subset of a Banach space $X$ and let $T: C \longrightarrow C$ be a self-mapping. Assume that $(x_n)_n$ is a sequence in $C$. $(x_n)_n$ is called an almost fixed point sequence (in short; a.f.p.s) for $T$ if
\vskip 0.3 cm
\centerline{$\displaystyle \lim_{n \longrightarrow + \infty}\|x_n - Tx_n \| = 0$.}
\end{definition}
\vskip 0.3 cm
\noindent It was proved (see Lemma 2.2 in \cite{Bet}) that if $C$ is bounded convex and $T$ is a Suzuki self-mapping on $C$ then $T$ has an a.f.p.s in $C$. In particular every nonexpansive mapping $T: C \longrightarrow C$ has an a.f.p.s.
\vskip 0.3 cm
\begin{definition}\label{def 3.3}\rm Let $C$ be a bounded closed convex subset of a Banach space $X$ and let $T: C \longrightarrow C$ be a self-mapping. $T$ is said to satisfy the condition $(L)$ if the following two conditions hold:
\begin{enumerate}
  \item If a subset $C_{0}\subset C$ is nonempty, closed, convex and $T$-invariant, then there exists an a.f.p.s for $T$ in $C_{0}$.
  \item For any a.f.p.s of $T$ in $C$ and all $x\in C$
\vskip 0.3 cm
\centerline{$\displaystyle \limsup_{n}\|x_{n}- Tx\|\leq \displaystyle \limsup_n \|x_{n}- x\|.$}
\vskip 0.3 cm
\end{enumerate}
\end{definition}

\begin{remark}\label{rem 3.6}\rm If $C$ is a weakly compact convex subset of a Banach space $X$. Then every generalized nonexpansive self-mapping on $C$ with $a+c > 0$ and every Suzuki self-mapping on $C$ satisfies the condition $(L)$. (see Propositions 3.4 and 3.6 in \cite{Fus}).
\end{remark}

\begin{remark}\label{rem 3.7}\rm The converse of Remark \ref{rem 3.6} is not true in general. Indeed, it was proved (see page 9 in \cite{Fus}) that the mapping $T$ of Example 2.6 satisfies condition $(L)$ but fails to be generalized nonexpansive or a Suzuki mapping.
\end{remark}

\noindent In the next result, we give a fixed point theorem concerning the class of mappings satisfying the condition $(L)$  which are nonexpansive with center $0$ in the setting of Banach spaces having KK property.

\begin{theorem}\label{th 3.4}\rm Let $X$ be a Banach space having KK property and let $C$ be a weakly compact convex subset of $X$. Assume that $T: C \longrightarrow C$ is a nonexpansive mapping with center 0 which satisfies the condition $(L)$. Then $T$ has a fixed point.
\end{theorem}

\proof If $0 \in C$, then the result is trivial and $0$ is a fixed point. Assume now that $0 \notin C$ and denote by $\Gamma$ the set $\{x \in C: \|x\|= \theta_{0}\}$ where $\theta_0 = \inf\{\|x\|: x \in C\} > 0$. So, $\Gamma$ is nonempty weakly compact convex subset of $X$. The fact that $\|Tx\|\leq \|x\|$ shows that $T(\Gamma)\subset \Gamma$. So, since $T$ satisfies the condition $(L)$ then $T$ has an a.f.p.s in $\Gamma$. Denote by $(x_{n})_{n}$ this a.f.p.s. But $\Gamma$ is weakly compact, thus from $(x_{n})_{n}$, we can extract a subsequence $(x_{n_{k}})_{k}$ in $\Gamma$ which converges weakly to some $y_{0}\in \Gamma$. Furthermore, from the definition of $\Gamma$,  for all integer $k$, we have $\|x_{n_{k}}\|= \|y_{0}\|= \theta_{0}> 0$. In addition, since $X$ satisfies KK property, we infer that $x_{n_{k}}$ converges in norm to $y_{0}$.
\noindent Now, by $(2)$ of the condition $(L)$ and using the fact that $(x_{n_k})_k$ is also an a.f.p.s for $T$, we get
\vskip 0.3 cm
\centerline{$0 \leq \displaystyle \liminf_{k} \|x_{n_k}- Ty_{0}\|\leq \displaystyle \limsup_{k}\|x_{n_k}- Ty_{0}\|\leq \displaystyle \limsup_{k}\|x_{n_{k}}- y_{0}\|= 0.$}
\vskip 0.3 cm
\noindent This leads to
\vskip 0.3 cm
\centerline{$\|y_{0}- Ty_{0}\|= \displaystyle \lim_{k}\|x_{n_{k}}- Ty_{0}\|= 0.$}
\vskip 0.3 cm
\noindent and so, $y_{0}= Ty_{0}$ which proves that $y_{0}$ is a fixed point for $T$ in $C$ and completes the proof.
\vskip 0.3 cm
\noindent Following Remark 3.6 and Theorem 3.4, we have

\begin{corollary}\label{cor 3.7}\rm Let $X$ be a Banach space having KK property and let $C$ be a weakly compact convex subset of $X$. If $T: C \longrightarrow C$ is nonexpansive with center $0$ satisfying one of the following assumptions:
\begin{description}
  \item[$(\imath)$] $T$ is a generalized nonexpansive mapping with $a+ c > 0$;
  \item[$(\imath\imath)$] $T$ is a Suzuki mapping.
\end{description}
\noindent Then $T$ has a fixed point in $C$.
\end{corollary}

\noindent Now, we are in position to state the following fixed point result.

\begin{theorem}\label{th 3.5}\rm Let $C$ be a weakly compact convex subset of a Banach space and let $T: C \longrightarrow C$ be a self-mapping. Then
\begin{description}
  \item[$(\imath)$] If $0 \in C$ and $T$ is a continuous asymptotically nonexpansive with center $0$. Then $0$ is a fixed point for $T$
  \item[$(\imath\imath)$] If $X$ has KK property and $0 \notin C$. Set $\theta_{0}= \inf\{\|x\|: x \in C\}$ and assume that $T$ satisfies the following assumptions:
\begin{description}
  \item[$(1')$] $T$ satisfies the condition $(L)$;
  \item[$(2')$] $T$ is asymptotically regular;
  \item[$(3')$] $T$ leaves the set $C_{\theta_{0}}= C \bigcap \{x \in X: \|x\|= \theta_{0}\}$ invariant.
\end{description}
\end{description}
\noindent Then $T$ has a fixed point in $C$.
\end{theorem}

\proof
\begin{description}
  \item[$(\imath)$] This claim is trivial (see the part $(\imath)$ in the proof of Theorem 3.2).
  \item[$(\imath\imath)$] Assume that $0 \notin C$. Since $C$ is weakly compact convex subset of $X$, then
\vskip 0.3 cm
\centerline{$\inf\{\|x\|: x \in C\}= \theta_{0} > 0$}
\vskip 0.3 cm
\noindent and the set $C_{\theta_{0}}$ is nonempty. On the other hand, by $(3')$ it is easy to observe that $C_{\theta_{0}}$ is a $T$-invariant weakly compact convex subset of $X$.
\noindent Let $x_{0} \in C_{\theta_{0}}$, since $C_{\theta_{0}}$ is weakly compact, from the sequence $(T^{n}x_{0})_{n}$, we can extract a subsequence $(T^{n_{k}} x_0)_{k}$ which converges weakly to some $y_{0} \in C_{\theta_{0}}$. Next, since $X$ has KK property, then $T^{n_{k}}x_{0} \longrightarrow y_{0}$ in norm. On the other hand, by $(2')$,  $T$ is asymptotically regular, then the sequence $(x_{n_k})$ defined by $x_{n_k} = T^{n_{k}} x_0$ is an a.f.p.s for $T$ in $C_{\theta_{0}}$. Now, from the fact that $T$ satisfies the condition $(L)$, we have
\vskip 0.3 cm
\centerline{$0 \leq \displaystyle \limsup_{k}\|x_{n_{k}}- Ty_{0}\|\leq \displaystyle \limsup_k \|x_{n_{k}}- y_{0}\|= 0.$}
\vskip 0.3 cm
\noindent Then $Ty_{0}= y_{0}$ which is the desired result.
\end{description}

\begin{remark}\label{rem 3.8}\rm Let $0 \notin C$. If $X$ satisfies KK property and assuming that $T: C \longrightarrow C$ is continuous and satisfies $(3')$ then $T$ has a fixed point in $C$. Indeed, by our proof, the set $C_{\theta_{0}}$ is a nonempty compact convex subset of $X$ and the result follows immediately from Schauder fixed point theorem. In this case, assumptions $(1')$ and $(2')$ of Theorem 3.5 can be dropped.
\end{remark}

\vskip 0.3 cm

\begin{remark}\label{rem 3.9}\rm It is worth noting that if
\vskip 0.3 cm
\centerline{$K= \{f \in L_{1}([0,1]): \ 0 \leq f \leq 2 \ a.e.  \displaystyle \int_{0}^{1}f(t)dt = 1\}$}
\vskip 0.3 cm
\noindent and $T$ is Alspach's mapping on $K$ given by
\vskip 0.3 cm
\centerline{$Tf(t)=\left\{
               \begin{array}{ll}
                 2f(2t) \wedge 2, \ \ \ & \hbox{if}\ 0\leq t \leq\displaystyle \frac{1}{2}, \\
                 (2f(2t-1)-2) \vee 0, \ \ \ & \hbox{if}\ \displaystyle \frac{1}{2}< t \leq1.
               \end{array}
             \right.$}
\vskip 0.3 cm
\noindent then $T$ is a free fixed point nonexpansive mapping which satisfies the assumption $(3')$ of Theorem 3.5 and in this case $\theta_{0}= 1$ with $K \bigcap C_{1}= K$. It was proved (see \cite{Als}) that $T$ is a nonlinear isometry satisfying $\|Tx\|= \|x\|$ for all $x \in K$. So, $T$ is nonexpansive with center 0. But, Banach space $L_{1}([0,1])$ does not have KK property.
\vskip 0.3 cm
\end{remark}

\noindent The following example shows that Theorem 3.5 is not true if we consider $C$, a bounded closed convex subset of an arbitrary Banach space even when assumptions $(1'), (2')$ and $(3')$ are satisfied.

\vskip 0.3 cm
\begin{example}\label{exa 3.4}\rm Let $X= C([0,1])$ equipped with the sup-norm. Let
\vskip 0.3 cm
\centerline{$C= \{x \in C([0,1]): 0=x(0)\leq x(t)\leq x(1)= 1\}$}
\vskip 0.3 cm
\noindent and let $T: C \longrightarrow C$ be the self-mapping defined by $Tx(t)= t x(t)$. Then $C$ is bounded closed and convex subset of $X$ and $T$ is a $(c)$-mapping (see Example in \cite{Bae}) (then it is asymptotically regular by Lemma 2.1). Obviously,
\vskip 0.3 cm
\centerline{$\|Tx\|\leq \displaystyle\sup_{t \in [0,1]}|tx(t)|\leq \displaystyle \sup_{t \in [0,1]}|x(t)|= \|x\|$}
\vskip 0.3 cm
\noindent and $T$ is nonexpansive with center $0$ which satisfies the condition $(3')$ with $\theta_{0}= 1$. However, $T$ is a free fixed point mapping. Notice that $X$ in this case does not have KK property.
\end{example}
\section{Applications}

\subsection{The case of alternate convexically nonexpansive mappings}
\noindent In \cite{Ami}, Amini-Harandi studied the existence of fixed points for a class of mappings called alternate convexically nonexpansive mappings defined on weakly compact convex subsets of strictly convex Banach spaces. To prove his result, Amini-Harandi used the existence of an afps (almost fixed point sequences) for such mappings and some properties of minimal sets associated with them (see p. 52 in \cite{Aks}). In \cite{Dow}, P. Dowling simplified Harandi's result without using any classical tool linked to the nonexpansive case and he observed that alternate convexically nonexpansive mappings form a subclass of that of nonexpansive mappings with center 0.
\vskip 0.3 cm
\begin{definition}\label{def 4.1}\rm Let $K$ be a nonempty subset of a Banach space $X$. A self-mapping $T: K \longrightarrow K$ is called alternate convexically nonexpansive if
\vskip 0.3 cm
\centerline{$\|\displaystyle \sum_{i= 1}^{n}\frac{(-1)^{i+1}}{n}Tx_{i}- Ty\| \leq \|\displaystyle \sum_{i=1}^{n}\frac{(-1)^{i+1}}{n}x_{i}- y\|$}
\vskip 0.3 cm
\noindent for all $n \in \mathbb{N}$ and $x_{i}, y \in K$.
\end{definition}

\noindent Let us state Amini-Harandi fixed point theorem (see \cite{Ami}).

\begin{theorem}\label{th 4.1}\rm Let $K$ be a weakly compact convex subset of a strictly convex Banach space and let $T: K \longrightarrow K$ be an alternate convexically nonexpansive mapping. Then $T$ has at least a fixed point.
\end{theorem}

\begin{remark}\label{rem 4.1}\rm To see that every alternate convexically nonexpansive mapping is nonexpansive with center 0, it suffices to take $n= 2, x_{1}, x_{2} \in K$ with $x_{1}= x_{2}$. So, the case of alternate convexically nonexpansive mappings becomes a particular case of the setting of nonexpansive mappings with center $0$. As a consequence, fixed point results associated with nonexpansive mappings with center $0$ hold also for alternate convexically nonexpansive mappings (see \cite{Dow, Fal2}).
\end{remark}
\vskip 0.3 cm
\noindent Furthermore, the same author in \cite{Dow} introduced the following weakening of the alternate convexically nonexpansiveness property.

\vskip 0.3 cm
\begin{definition}\label{def 4.2}\rm Let $K$ be a nonempty subset of a Banach space $X$. A self-mapping $T: K \longrightarrow K$ is called $k$-alternate convexically nonexpansive if
\vskip 0.3 cm
\centerline{$\|\displaystyle \sum_{i= 1}^{n}\frac{(-1)^{i+1}}{n}Tx_{i}- Ty\| \leq \|\displaystyle \sum_{i=1}^{n}\frac{(-1)^{i+1}}{n}x_{i}- y\|$}
\vskip 0.3 cm
\noindent for all $1 \leq n \leq k$ and $x_{i}, y \in K$.
\end{definition}
\vskip 0.3 cm
\noindent It was observed (see Remark 3 in \cite{Dow}) that if $T$ is $k$-alternate convexically nonexpansive then $T$ is nonexpansive with center 0. In addition, in the same paper it was proved that Alspach transformation (see Remark 3.9) is an example of 1-alternate convexically nonexpansive that is not 2-alternate convexically nonexpansive.
\vskip 0.3 cm
\noindent From Corollary 3.2, we can establish the following fixed point result concerning 2-alternate convexically nonexpansive mappings.
\vskip 0.3 cm
\begin{corollary}\label{co 4.1}\rm Let $K$ be a closed convex subset of a reflexive strictly convex Banach space and let $T: K \longrightarrow K$ be a 2-alternate convexically nonexpansive mapping. Then $T$ has at least a fixed point.
\end{corollary}

\vskip 0.3 cm
\begin{remark}\label{re 4.2}\rm Corollary 4.1 extend Theorem 3 of \cite{Dow} to the case of unbounded closed convex subsets of reflexive strictly convex Banach spaces.
\end{remark}

\vskip 0.3 cm
\subsection{The linear isometries $l_g$ and $r_g$}
\noindent We start this section by investigating some exotic situations associated with linear isometries acting on $L_{p}(G), 1< p < \infty$ where $G$ is a locally compact group.

\vskip 0.3 cm
\noindent Let $G$ be a locally compact group with a left Haar measure $v$ and modular function $\Delta$ defined by
\vskip 0.3 cm
\centerline{$\Delta(g) \displaystyle \int_{G}k(xg)dv(x)= \displaystyle \int_{G}k(x)dv(x)$}
\vskip 0.3 cm
\noindent for $k \in C_{0}^{c}(G)$ the space of continuous functions $k$ vanishing off compact subsets of $G$. The left and the right translations in $L_{p}(G), 1< p < \infty$ by $g \in G$ are given respectively by
\vskip 0.3 cm
\centerline{$l_{g}f(x)= f(gx) \ \hbox{and} \ (r_{g}f)(x)= \Delta^{\frac{1}{p}}(g)f(xg)$}
\vskip 0.3 cm
\noindent for all $x \in G$. These mappings satisfy $l_{g_{1}}l_{g_{2}}= l_{g_{2}g_{1}}$ and $r_{g_{1}}r_{g_{2}}= r_{g_{1}g_{2}}$ for all $g_{1}, g_{2}\in G$. Furthermore each $l_{g}$ and $r_{g}$ is a linear isometry.
\vskip 0.3 cm
\noindent A subset $C_0 \subset L_p(G)$ is called left (resp. right) invariant if $l_a(C_0) \subset C_0$ (resp. $r_a(C_0) \subset C_0)$ for each $a \in G$.
\vskip 0.3 cm
\noindent In \cite{Lau}, A. T-M. Lau studied closed convex left or right invariant subsets of $L_p(G)$. He proved in particular that if $G$ is a locally compact noncompact group then every closed convex left invariant subset $C_0$ of $L_p(G)$ must contain 0. In addition, if $C_0$ is assumed to be compact convex, then $C_0$ is reduced to the singleton $\{0\}$.

\vskip 0.3 cm
\noindent If $K$ is a nonempty closed convex subset of $L_p(G) (1 < p < \infty)$ which is invariant by every $l_g$ (resp. $r_g$) ($g \in G$),  we denote by $F(l_{g}), g\in G$ (resp. $F(r_{g}), g\in G$)  the set of fixed points of $l_{g}$ (resp. $r_g$) in $K$.
\vskip 0.3 cm
\noindent First of all, we remark that the family of mappings $(l_{g}, g \in G)$ (resp. $r_{g}, g \in G)$ are not commuting in general. But since for $1 < p < \infty, L_{p}(G)$ is strictly convex (see page 293, Corollary 20.14 in \cite{Hew}) and the fact that $\|l_{g}(x)\|= \|x\|$ for all $x \in L_{p}(G)$ and $g \in G$ then we deduce that if $K$ is a closed convex invariant subset of $L_{p}(G) (1 < p < \infty)$ then each $l_g$ ($ g \in G)$ is nonexpansive with center 0. By using Corollary 3.2, we have $l_{g}(x_{0})= x_{0}$ for all $g\in G$,  so if
\vskip 0.3 cm
\centerline{$x_{0}\in \displaystyle \bigcap_{g \in G}F(l_{g})$}
\vskip 0.3 cm
\noindent and $x_{0}\neq 0$, then by Lau's result indicated above,  necessarily $G$ is a compact group which is a contradiction. Thus,  we can derive the following.
\vskip 0.3 cm
\begin{corollary}\label{cor 4.2}\rm Let $G$ be a locally compact group. Assume that $G$ is noncompact and there exists a closed convex subset of $L_{p}(G), 1 < p < \infty$ which is invariant by each $l_{g}, g \in G$. Then
\vskip 0.3 cm
\centerline{$\displaystyle \bigcap_{g \in G} F(l_{g})= \{0\}.$}
\end{corollary}
\vskip 0.3 cm
\noindent In the next result, we give another characterization of the noncompactness of a locally compact group $G$ by means of orbits associated with the mappings $l_g$ or $r_g$. We will restrict our proof to the case of $l_g (g \in G)$ mappings.
\vskip 0.3 cm

\begin{corollary}\label{cor 4.3}\rm Let $G$ be a locally compact group. Assume that $K$ is an arbitrary nonempty weakly compact convex subset of $L_{p}(G), 1 < p< \infty$ which is $l_g$-invariant for all $g \in G$. Then the following assertions are equivalent:
\begin{description}
  \item[$(\imath)$] $G$ is noncompact;
  \item[$(\imath\imath)$] For all fixed $h \in K$, we have $0 \in \overline{co}\{l_{{g}} h : g \in G\}$.
\end{description}
\end{corollary}

\proof \begin{description}
         \item[$(\imath)\Longrightarrow (\imath\imath)$] Assume that there exists $h_0 \in K$ such that $0 \notin \overline{co}\{l_{g} h_0: g \in G\}$. Since each $l_{a}(a \in G)$ is continuous and affine (since it is linear), $\overline{co}\{l_{g}h_0: g \in G \}$ is a closed convex subset which is invariant by each $l_{a}, a \in G$. This fact contradicts Lau's result.
         \item[$(\imath\imath)\Longrightarrow (\imath)$] Assume that $G$ is compact, then if we take $f_{0} =1$ (that is $f_0(x)= 1, \forall x \in G$), then $f \in L_{p}(G), 1 < p< \infty$, we have $\overline{co}\{l_{g}f_{0}: g \in G\}= \{f_{0}\}$ which does not contain the origin. By taking $K = \{f_0\}$,  this is a contradiction.
       \end{description}
       \vskip 0.3 cm

\begin{remark}\label{rem 4.3}\rm Recall that the paper \cite{Lim} is an important investigation on the existence of fixed points for isometries defined on weakly compact convex subsets of Banach spaces. Indeed, in the indicated paper, the authors proved that isometries which are defined on bounded closed convex subsets of uniformly convex Banach spaces have the Chebychev center as a common fixed point. It is easily seen that in our setting related to the isometries $l_g$ or $r_g$ ($g \in G)$ this Chebychev center is reduced to the set $\{0\}$.

\end{remark}

\subsection{The case of a nonlinear transport equation}
\noindent Here, using our results in Section 1, we will investigate the existence of a solution for the following boundary problem
\begin{align}
\lambda\varphi(x,v)+ v. \nabla_{x}\varphi(x,v)+ \sigma(v)\varphi(x,v)= &  \int_{V}h(x,v,v')f(x,v',\varphi(x,v'))d\mu(v')  \\
                 \varphi_{-} =  & H(\varphi_{+}) \nonumber
\end{align}

\noindent where $\lambda \in \mathbb{R}$, $f(.,.,.)$ is a measurable nonlinear function of $\varphi$ and $h(.,.,.)$ is a measurable function from $D\times V \times V$ to $\mathbb{R}$ where $D$ is a smooth open subset of $\mathbb{R}^{n}$ that represents the domain of positions and $V$ is the support of the Radon measure $\mu$ on ${\mathbb{R}}^{n}$ with $\mu(\{0\})= 0$. Recall that $V$ is the velocities space. The unknown function $\varphi(x, v)$ is the number (or probability) density of gas particles having the position $x$ and the velocity $v$. The homogeneous function $\sigma(.)$ is called the collision frequency. The boundary conditions are modeled by
\vskip 0.3 cm
\centerline{$\varphi_{-}= H (\varphi_{+})$}
\vskip 0.3 cm
\noindent where $\varphi_{-}$ (resp. $\varphi_{+}$) is the restriction of $\varphi$ to $\Gamma_{-}$ (resp. $\Gamma_{+}$) which is the incoming (resp. outcoming) part of the phase space boundary and $H$ is a bounded linear operator acting between suitable Lebesgue function spaces on $\Gamma_{+}$ and $\Gamma_{-}$, covering in particular the classical boundary conditions (vacuum boundary conditions corresponding to $H= 0$, periodic boundary conditions, reflexive boundary conditions,...)
\vskip 0.3 cm
\noindent In our setting, the function $h(.,.,.)$ is chosen such that the linear operator
\begin{align}
R: L_p(D \times V) \longrightarrow & L_{p}(D \times V)(1< p < \infty) \nonumber \\
\varphi \longrightarrow & \int_{V}h(x,v,v')\varphi(x,v') d\mu(v') \nonumber
\end{align}
\noindent is bounded.

\begin{definition}\label{def 4.3}\rm A function $g: D\times V\times \mathbb{R} \longrightarrow \mathbb{R}$ is a Caratheodory function if the following condition is satisfied
\begin{align}
\hbox{for all} \ \ s \in \mathbb{R},  (t,s) \longrightarrow & g(t,s) \ \ \hbox{is measurable in} \ \   D\times V \nonumber \\
s \longrightarrow & g(t,s) \ \ \hbox{is continuous on} \ \mathbb{R} \ a.e. \ t\in D \times V. \nonumber
\end{align}
\end{definition}

\begin{remark}\label{def 4.4}\rm If $f$ is a Caratheodory function, then we can define the Nemytskii operator $N_{f}$ by
\vskip 0.3 cm
\centerline{$(N_{f}\varphi)(x,v)= f(x,v, \varphi(x,v))$}
\vskip 0.3 cm
\noindent for all $(x,v)\in D \times V$. In addition, if the operator $N_{f}$ acts on $L_{p}(D \times V)$, then $N_{f}$ is continuous and it takes bounded sets into bounded sets.
\end{remark}

\begin{remark}\label{def 4.5}\rm It is easily seen that if $N_{f}$ is a nonexpansive mapping with center $0$ on $L_{p}(D \times V)$ then $N_{f}$ takes every ball $\overline{B}(0, r)$ in $L_{p}(D \times V)$ into itself (and consequently has 0 as a fixed point in $\overline{B}(0, r)$).
\end{remark}

\noindent In \cite{Lat}, the author studied the existence of solutions for the nonlinear equation $(4.1)$ by using some compactness results in transport theory which require the boundedness and the convexity of $D$ together with the regularity of the bounded linear operator $R$ (that is the compactness on $L^{p}(V)$, if the position $x$ is fixed). So, the solution is derived from Schauder's fixed point theorem for convenable mappings acting on balls with center $0$.
\vskip 0.3 cm
\noindent In our main results below, compactness assumptions are not required.
\vskip 0.3 cm
\noindent Denote by $T_{H}$ the unbounded linear operator defined on $L_{p}(D \times V), 1 < p < \infty$ by
\begin{align}
    T_{H}\varphi(x,v)=&  - v \nabla_{x}\varphi(x,v)- \sigma(v)\varphi(x,v), \nonumber \\
\varphi_{-}= & H(\varphi_{+}). \nonumber
\end{align}

\noindent Our assumptions denoted by $({\cal H})$ are the following:
\begin{itemize}
  \item $D$ is an open smooth subset of $\mathbb{R}^{n}$.
  \item $R$ is a bounded operator on $L_{p}(D \times V)$.
  \item $f$ is a Caratheodory function.
  \item $N_{f}$ acts from $L_{p}(D \times V)$ into itself $(1< p <\infty)$.
  \item There exists $r_{0}$ and $x_{0}\in L_{p}(D \times V)$ such that $N_{f}$ is a nonexpansive self-mapping with center $0$ on $\overline{B}(x_{0}, r_{0})$.
  \item For $\lambda$ sufficiently large and $\|H\| < 1$, the mapping $B_{\lambda}= (\lambda - T_{H})^{-1}R$ leaves $\overline{B}(x_{0}, r_{0})$ invariant.
\end{itemize}

\begin{theorem}\label{th 4.2}\rm Assume that $({\cal H})$ is satisfied. Then there exists $\lambda_{0} > 0$ such that for all $\lambda > \lambda_{0}$, the problem $(4.1)$ has at least one solution in $\overline{B}(x_{0}, r_{0})$.
\end{theorem}

\proof It is easy to observe that the problem $(4.1)$ has a solution if and only if the mapping
\vskip 0.3 cm
\centerline{$S_{\lambda}= (\lambda- T_{H})^{-1}RN_{f}$}
\vskip 0.3 cm
\noindent has a fixed point. By our assumption $S_{\lambda}$ leaves $\overline{B}(x_{0}, r_{0})$ invariant. On the other hand, for all $\varphi \in \overline{B}(x_{0}, r_{0})$,  since $N_f$ is nonexpansive with center 0, we have
\vskip 0.3 cm
\centerline{$\|S_{\lambda}\varphi\|\leq \|(\lambda- T_{H})^{-1}RN_{f}\varphi\|\leq \|(\lambda- T_{H})^{-1}\|\|R\|\|N_{f}\varphi\|\leq \|(\lambda- T_{H})^{-1}\|\|R\| \|\varphi\|$}
\vskip 0.3 cm
\noindent Following Lemma 2.2 in \cite{Lat}, we have
\vskip 0.3 cm
\centerline{$\|(\lambda- T_{H})^{-1}\|\leq \displaystyle \frac{1}{\lambda+ \lambda^{*}}$}
\vskip 0.3 cm
\noindent when $\lambda^{*}= - \displaystyle \lim_{|v|\longrightarrow 0} \sigma(v).$
\vskip 0.3 cm
\noindent Now, since
\vskip 0.3 cm
\centerline{$\displaystyle \lim_{\lambda \longrightarrow + \infty}\frac{ \|R\|}{\lambda+ \lambda^{*}}=0$,}
\vskip 0.3 cm
 \noindent there exists $\lambda_{0} \in \mathbb{R}$ such that for all $\lambda \in \mathbb{R}$ satisfying  $\lambda \geq \lambda_0$, we infer that
\vskip 0.3 cm
\centerline{$\displaystyle \frac{ \|R\|}{\lambda+ \lambda^{*}}\leq 1$}
\vskip 0.3 cm
\noindent which gives that
\vskip 0.3 cm
\centerline{$\|S_{\lambda}\varphi\|\leq \|\varphi\|$}
\vskip 0.3 cm
\noindent for all $\varphi \in \overline{B}(x_{0}, r_{0})$ and $\lambda \in \mathbb{R}$ satisfying $\lambda \geq \lambda_0$.
\noindent So, $S_{\lambda}$ is a nonexpansive mapping with center $0$. In addition, since for $ 1< p< + \infty$, Banach spaces $L_{p}(D \times V)$ are reflexive and strictly convex, then $\overline{B}(x_{0}, r_{0})$ is a closed convex subset of $L_p(D \times V)$ and the fixed point for $S_{\lambda}$ follows from Corollary 3.2.

\begin{remark}\label{rem 4.6}\rm Theorem 4.2 does not require any compactness results or a specific geometry on the spaces of positions $D$.  In addition, following Corollary 3.2, in conditions $(5)$ and $(6)$ of $({\cal H})$, we can replace the closed ball $\overline{B}(x_0, r)$ by any nonempty closed convex subset.
\end{remark}

\begin{remark}\label{rem 4.7}\rm In the case where $D= \mathbb{R}^{3}$ and $V= \mathbb{R}^{3}$ then $D \times V = \mathbb{R}^{3}\times \mathbb{R}^{3}$ is a locally compact group. When $\sigma(.)\equiv 0$, then if $H = 0$, $T_{H}$ generates a semigroup of contractions on $L_{p}(\mathbb{R}^{3}\times \mathbb{R}^{3}),1 < p < +\infty $ given by
\vskip 0.3 cm
\centerline{$U(t)\varphi(x,v)= \varphi(x-tv,v) \hskip 1 cm (x,v)\in \mathbb{R}^{3}\times \mathbb{R}^{3}.$}
\vskip 0.3 cm
\noindent It is interesting to investigate weakly compact (resp. compact) convex subsets which are invariant by the flow $(U(t))_{t \geq 0}$ on $L_{p}(\mathbb{R}^{3}\times \mathbb{R}^{3})$.
\end{remark}

\noindent {\bf \it Funding:} This work is supported by the research team RPC (Controllability and Perturbation Results) in the laboratory of Informatics and Mathematics (LIM) at the university of Souk-Ahras (Algeria).

\end{document}